\documentclass[11pt]{amsart}
\usepackage{geometry}                
\usepackage[parfill]{parskip}    
\usepackage{graphicx}
\usepackage{amssymb}
\usepackage{epstopdf}
\usepackage{xcolor}

\DeclareFontFamily{U}{wncy}{}
    \DeclareFontShape{U}{wncy}{m}{n}{<->wncyr10}{}
    \DeclareSymbolFont{mcy}{U}{wncy}{m}{n}
    \DeclareMathSymbol{\Sha}{\mathord}{mcy}{"58} 

 \newcommand{\R}{\mathbb R}
 \newcommand{\Z}{\mathbb Z}
  \newcommand{\G}{\mathfrak G }

 \newcommand{\beq}{\begin{equation}}
 \newcommand{\eeq}{\end{equation}}
\newcommand{\beqq}{\begin{equation*}}
 \newcommand{\eeqq}{\end{equation*}}

\newcommand{\M}{{\text{M}}}
\newcommand{\NN}{\mathbb N}
\newcommand{\CC}{\mathbb{C}\,}

\newcommand{\RR}{\mathbb R}
\newcommand{\ZZ}{\mathbb Z}

\newcommand{\Lt}{{L^2(\RR)}}
 \newcommand{\PW}{PW_1}
  \newcommand{\PWa}{PW_{1/\alpha}}
 \newcommand{\PWb}{PW_a}

\newtheorem{Thm}{Theorem}[section]
\newtheorem{theorem}[Thm]{Theorem}

\newtheorem{corollary}[Thm]{Corollary}
\newtheorem{remark}[Thm]{Remark}

\newtheorem{definition}[Thm]{Definition}

\title{Gabor frame operator for the Cauchy kernel}
\author{Yurii Belov}
\address{Yurii Sergeevich Belov,
\newline Department of Mathematics and Computer Sciences, 
St.~Petersburg State University,  29 Line 14th (Vasilyevsky Island), 199178, St. Petersburg, Russia,
\newline {\tt j\_b\_juri\_belov@mail.ru} }
\author{Aleksei Kulikov}
\address{Aleksei Igorevich Kulikov, 
\newline Tel Aviv University, School of Mathematical Sciences, Tel Aviv, 69978, Israel,
\newline {\tt lyosha.kulikov@mail.ru} }
\author{Yurii Lyubarskii}
\address{Yurii Ilích  Lyubarskii,  
 \newline Department of Mathematical Sciences, 
Norwegian University of Science and Technology, NO-7491 Trondheim, Norway,}
\smallskip
\address{Institute of Mathematics, Ufa Federal Research Center, RAS,
Chernyshevky str. 112,
450008, Ufa, Russia
\newline {\tt yura@math.ntnu.no} }

\thanks{\sc Yu. Belov, A.Kulikov, Yu. Lyubarskii, Gabor frame operator for the Cauchy kernel.}
\thanks{\copyright \ Belov Yu.,  Kulikov A., Lyubarskii Yu. \ 2022}
\thanks{\rm The research of Yu. Lyubarskii was supported by the grant of Russian Science Foundation (project no.  21-11-00168)}
\thanks{\it Submitted August 15, 2022.}
\begin{document}
\maketitle
\begin{abstract}
 
We obtain frame bounds estimates and 
 the Gabor frame operator $S=S^{\alpha,\beta}$ for Gabor frames generated by the Cauchy kernel. 
 In addition we find the explicit expression   for the canonical dual window for all values of the
  lattice parameters $\alpha,\beta$, $\alpha\beta\leq 1$. 

\medskip



\end{abstract}
\section{Introduction}

In this  paper we continue the study {of} the Gabor frames generated by the Cauchy kernel
 $$  
  g(t):=
  \frac1 {t-iw},  \  w\in \CC\setminus  i \RR.
  $$
 In \cite{BKL} we described all irregular lattices $\Lambda \times \M$, $\Lambda \in \RR$, $\M\in \RR$ such 
 that the corresponding collection of the time-frequency shifts 
 $\{e^{2\pi i\mu t}g(t-\lambda)\}_{\lambda\in  \Lambda, \mu\in \M}$ {is} a frame in $L^2(\RR)$.
 
 For the regular lattices $\alpha \ZZ \times \beta \ZZ$ {we} 
  obtain more detailed 
 information, namely the explicit expression of the frame operator and the canonical dual window as well as estimates {for} the frame bounds. To the best of our knowledge explicit expressions for the canonical dual windows were so far known only
 for the Gabor frames   which  are Riesz  bases.
 
\subsection*{Gabor frame bounds} We {recall} the main notation. Given $\alpha, \beta > 0$ we denote by $g_{m, n}$ {\it time-frequency shift } of $g$, i.e.
$$
g_{m,n}(t)= 
\frac{ e^{2i\pi \beta n t}}{t-\alpha m - iw}, \ m,n \in \Z,
$$
and consider the corresponding Gabor system
$$
\G(g; {\alpha, \beta}) =\left \{g_{m,n}\right \}_{m,n \in \Z}.
$$
 
Taking if need be the complex conjugation and then shift in $t$ we {can} always assume {that $w\in \R$ and} $w>0$. It is known (see e.g. \cite{Ja, GRS,BKL})  that $\G(g; {\alpha, \beta}) $ forms a frame in $\Lt$  for all
$\alpha \beta\leq 1$, that is
there exist constants $0< A(\alpha, \beta; w)\leq B(\alpha, \beta;w)<\infty$, such that  
\beq
\label{eqno01}
A(\alpha, \beta; w) \| f\|^2  \leq \sum_{m,n\in \Z} | \langle f, g_{m,n} \rangle |^2 \leq   B(\alpha, \beta; w) \| f\|^2 ,
\eeq
for all $f\in \Lt $. These constants are responsible for stability of the
reconstruction of  functions  in $ \Lt$ through the {\em frame coefficients}  
$\langle f,g_{m,n} \rangle$.
They are called  {\em lower} and {\em upper} {\em frame constants} respectively.   From now on we write $A(\alpha,\beta;w)$, $B(\alpha,\beta;w)$ as the largest $A$ and the smallest $B$ respectively such that (\ref{eqno01}) holds for all $f\in L^2(\mathbb{R})$.

Relation (\ref{eqno01}) can be expressed in terms of the {\it frame operator} $S^{\alpha,\beta}_g: \Lt \to \Lt$,
\beq
\label{eqno02}
 S^{\alpha,\beta}_g f(t)=\sum_{m,n\in \ZZ} \langle f,g_{m,n} \rangle g_{m,n}(t), \ f\in \Lt.
 \eeq
 We should have 
 $$
 A(\alpha, \beta; w)Id \le S^{\alpha,\beta}_g \le  B(\alpha, \beta;w)Id.
 $$
 
 In  the special case of integer oversampling $(\alpha \beta)^{-1} \in \NN$  these constants have been  
 studied in \cite{Ja, FS}. 
In particular, on the critical hyperbola $\alpha\beta=1$ the {\em optimal} values of $A$ and $B$ {can} be found 
  by using the  {\it Zak transform}  of the window $g$. 
We refer the reader to \cite{Gro}[Ch. 8]
for  the basic facts on the Zak transform. We remind that Zak transform is defined by the relation 
$$
\mathcal{Z}_\alpha : f \mapsto \mathcal{Z}_\alpha f(t, \omega) = \sum_{k=-\infty}^\infty f(t-\alpha k)e^{2\pi i \alpha k \omega}, \ f\in \Lt.
$$ 
This is a unitary mapping of $\Lt$ onto $L^2((0,\alpha)\times (0,1/\alpha))$, and, for $\alpha\beta=1$, the frame operator $S_g$ is unitary equivalent to the multiplication by $\alpha ^{1/2}\mathcal{Z}_\alpha g$ in  $L^2((0,\alpha)\times (0,1/\alpha))$. Therefore, for $\alpha\beta=1$, we obtain the optimal values of the frame constants. 
$$
 A(\alpha, 1/\alpha; w) =\alpha\inf_{t\in [0,\alpha],\omega\in [0,1/\alpha]} |\mathcal{Z}_g(t,\omega)|^2,\quad B(\alpha, 1/\alpha;w)  =\alpha\sup_ {t\in [0,\alpha],\omega\in [0,1/\alpha]} |\mathcal{Z}_g(t, \omega)|^2.
 $$ 
In the case of the Cauchy kernel the Zak transform can be found explicitly.  We have (after regularization)
$$
\mathcal{Z}_\alpha g(t,\omega)=\frac{{-2\pi i}}{\alpha}e^{2\pi i t\omega}\frac{e^{2\pi \omega {w}}}{1-e^{2\pi i t\slash\alpha}e^{2\pi w\slash\alpha}}, \quad t\in[0,\alpha],  \omega\in[0,\alpha^{-1}].
$$

Hence, for $\alpha\beta=1$,  
$$
A(\alpha, 1/\alpha; w)  = \frac{ 4 \pi^2}{\alpha}\biggl{(}\frac{1}{e^{2\pi w\slash \alpha}+1}\biggr{)}^2,
\quad B(\alpha, 1/\alpha; w)  = \frac{4 \pi^2}{\alpha}\biggl{(}\frac{e^{2\pi w\slash \alpha}}{e^{2\pi w/\alpha}-1}\biggr{)}^2.
$$


\begin{theorem}
 For all $\alpha, \beta > 0$, $\alpha \beta \leq 1$ the   frame constants,  $A(\alpha, \beta;w)$,  $B(\alpha, \beta;w)$ admit the estimates 
\beq 
\label{fi}
A(\alpha, \beta;w) \geq  \frac{4\pi^2 }{\alpha} \left (\frac {1}{e^{2\pi \beta w}+1}\right )^2,  \
B(\alpha, \beta;w) \leq  \frac{4\pi^2 }{\alpha} \left (\frac {e^{2\pi \beta w}}{e^{2\pi \beta w}-1}\right )^2
 \eeq
\label{main}
\end{theorem}

Proof of Theorem \ref{main} involves study of the corresponding T{oe}plitz operators as well as  interpolation in spaces of entire functions.

\subsection*{Frame operator and  the canonical dual window}

  For any Gabor frame $\G(g; {\alpha, \beta})  =\left \{g_{n,m}\right \}_{n,m \in \Z}$ there exists {a} {\it dual window } $h\in L^2(\mathbb{R})$  (see e.g. \cite{Gro}[Ch. 5]) such that any $f\in \Lt$ admits the representation   
 $$
 f=\sum_{m,n}\langle f,g_{m,n} \rangle h_{m,n}=\sum_{m,n}\langle f,h_{m,n} \rangle g_{m,n}.
 $$
 The dual window is by no mean unique (unless we are dealing with a Riesz basis), yet there is {\em the canonical dual window} $\gamma$ which minimises  the $l^2$ norm of the coefficients $\{\langle f,h_{m,n} \rangle \}$ such that the representation 
 $$
 f=\sum_{m,n}\langle f,g_{m,n} \rangle \gamma_{m,n}=\sum_{m,n}\langle f,\gamma_{m,n} \rangle g_{m,n}
 $$
becomes in a sense optimal.
 
 The canonical dual window is defined by the relation
 $$
 \gamma = S^{-1} g,
 $$
 so far its precise expression was unknown for the Gabor frames with the density above the critical one.  
 
 In the case of the Cauchy kernel we will find the frame operator $S$, its inverse $S^{-1}$, and the explicit expression of the canonical dual window for all values $\alpha$, $\beta$, $\alpha\beta\leq 1$

We will use standard form of the Fourier transform
$$
\hat{f}(\xi)= \int_{-\infty}^\infty f(t)e^{-2i\pi \xi t} dt.
$$

Let us consider distriburion $h=h_{\alpha,\beta}$ such that 

$$\hat{h}(\xi)={\pi}\sum_m\chi_{[-(2\alpha)^{-1},(2\alpha)^{-1}]+\beta m}(\xi)e^{4\pi w (\xi-\beta m)},$$
{that is
$$h(x) = \frac{\sin(\frac{\pi}{\alpha}(x-2iw))}{x-2iw}\Sha(\beta x),$$
where $\Sha(x)$ is the Dirac comb $\Sha(x) = \sum_{n\in \Z} \delta_n$.}

\begin{theorem}  The frame  operator $S=S^{\alpha,\beta}$  for Cauchy kernel $g$ is given by 
$$
Sf(x) =\frac{\pi}{\alpha\sin\left(\frac{\pi}{\alpha}(x-iw)\right)}\left(
\frac{f(t)}{\sin\left(\frac{\pi}{\alpha}(t+iw)\right)}\ast h\right).
$$
 The inverse   operator $S^{-1}$ is given by
$$S^{-1}f(t)=\frac{\alpha}{\pi}\sin\left(\frac{\pi}{ \alpha}(t+iw)\right)\left(f(x)\sin\left(\frac{\pi}{ \alpha}(x-iw)\right)\ast \tilde{h}\right),$$
where $\tilde{h}$ is such that $\hat{\tilde{h}}\hat{h}=1$.

 In particular, for $\alpha\beta\geq 1\slash 2$ we have
$$\gamma(t)=\frac{\alpha}{\pi} \sin\left(\frac{\pi}{\alpha}(t+iw)\right)\biggl{[}
\frac{\sin(\pi(\beta-\varepsilon)(t{+}iw))}{{\pi(t+iw)}}+ 
\frac{\cos(\pi \beta(t{-}iw))}{\cosh(2\pi w\beta)}
\frac{\sin(\pi\varepsilon(t{+}iw))}{{\pi(t+iw)}}\biggr{]}$$
where $\varepsilon=\frac{1}{\alpha}-\beta$.
\label{main2}
\end{theorem}

Here we understand convolution of $h$ and $\tilde{h}$  with the functions in $\Lt$ in the distributional sense  i.e. as the multiplication  by $\hat{h}$ and $\hat{\tilde{h}}$  respectively  in the Fourier domain.
 For the functions  $f\in C_0^\infty(\RR)$ this of course becomes the regular convolution.

{We are able to compute the function $\gamma$ for all pairs $\alpha, \beta$ with $\alpha\beta \le 1$ but for the sake of readers's sanity we decided not to present it in the introduction. See \eqref{even n},  \eqref{odd n} for the formula for the function $\gamma$ in the general case.} 

Theorem \ref{main2}   also allows us to { establish new} estimates for frame constants.

\begin{corollary} We have
 $$
 \frac{4\pi}{\alpha}\frac{\inf_\xi|\hat{h}_\beta(\xi)|}{(e^{\pi w\slash \alpha}+e^{-\pi w\slash \alpha})^2}\leq A(\alpha, \beta;w)\leq \frac{4\pi}{\alpha}\frac{\inf_\xi|\hat{h}_\beta(\xi)|}{(e^{\pi w\slash \alpha}-e^{-\pi w\slash \alpha})^2}, 
 $$

$$
 \frac{4\pi}{\alpha}\frac{\sup_\xi|\hat{h}_\beta(\xi)|}{(e^{\pi w\slash \alpha}+e^{-\pi w\slash \alpha})^2}\leq B(\alpha, \beta;w)\leq \frac{4\pi}{\alpha}\frac{\sup_\xi|\hat{h}_\beta(\xi)|}{(e^{\pi w\slash \alpha}-e^{-\pi w\slash \alpha})^2}.
 $$
 
\label{cor}
\end{corollary}

These estimates   in some cases better 
 than ones in Theorem \ref{main}. For critical hyperbola $\alpha\beta=1$  the lower estimate for $A$ and upper estimate for $B$ are also sharp.

It is worth mentioning that the Fourier transform of $\gamma$  is supported on the interval $[-\frac{1}{\alpha}, \frac{1}{\alpha}]$ (for all $\alpha, \beta$, $\alpha\beta\leq1$) and, moreover, $\gamma(t) = \frac{Q(t)}{t+iw}$ for some trigonometric polynomial $Q$  vanishing at $-iw$.   We do not have an a priori explanation why this must be true. In particular, this means that dual window for one-sided exponential $e^{-2\pi w \xi}\chi_{>0}(\xi)$, $w>0$ (Fourier transform of the Cauchy kernel) has compact support. 
  The size of this support depends upon the lattice parameters, of course. 

\section{Frame coefficients}
\subsection*{a} { Normalization.} Recall that we are already dealing  with  the case $w>0$. A direct calculation shows  
\beq 
\label{norm}
A(\alpha, \beta;w)= \beta A(\alpha\beta, 1, \beta w), B(\alpha, \beta)= \beta B(\alpha\beta, 1; \beta w).
\eeq
 Therefore  we  start with  study $A(\alpha, 1;w)$, and  $B(\alpha, 1;w)$,  $0<\alpha \leq 1 $ only and 
 later formulate   the general result.  
\medskip
 
 \begin{remark}  For the complex values of $w\in \CC\setminus i \RR $,    one has to replaced $w$ by 
 $|\Re w|$ { in the above relations}.
\end{remark}

 \medskip
\subsection*{b} {Calculation of the frame coefficients. } 
We need the following

\medskip

\begin{definition}
Given $a>0$,  the Paley-Wiener space $\PWb$ is defined as
\beqq
\label{3}
\PWb=\left \{f:\, f(z)=\int_0^{a}e^{2i\pi\xi z} \hat{f}(\xi)d\xi, \, \hat{f}\in L^2(0, a)\right \}.
\eeqq 
\end{definition}

This space consists of entire functions of exponential type, which belong to $\Lt$ and have the indicator diagram included in $[-2ia \pi,0]$.
We refer the reader to \cite{L} for the detailed description of the Paley-Wiener spaces as well as for other facts on entire functions. 

We fix now $\alpha\in (0,1]
$, so  
$$
g_{m,n}(t)=
 \frac{e^{2i\pi n t}}{t-\alpha m - iw}, \ m,n \in \Z.
$$
 Given $f\in \Lt$,  we set  
 
 \beq
 \label{5}
 c_{m,n}=\frac 1{2i\pi}\langle f, g_{m,n}\rangle = \frac 1{2i\pi}\int_{-\infty}^\infty f(t)
       \frac{e^{-2i\pi n t}}{t-\alpha m{+}iw}dt.
 \eeq
 The frame inequality now takes the form  
 $$
   A(\alpha, 1; w) \| f\|^2  \leq 4\pi^2 \sum_{m,n\in \Z} | c_{m,n}  |^2 \leq     B(\alpha, 1; w) \| f\|^2 ,
$$

Let
 \beqq
 \label{6}
 f_k(t)=\int_{k}^{{k+1}}\hat{f}(\xi)e^{2i\pi\xi t}d\xi.
 \eeqq

We then have 
 \beq
 \label{7}
 f(t)=\sum_{k=-\infty}^\infty f_k(t), \ \|f\|^2=\sum_{k=-\infty}^\infty\|f_k\|^2.
 \eeq
  
Denote  also 
\beqq
 h_k(t):= e^{-2i\pi kt}f_k(t)\in \PW. 
\eeqq

 A straightforward calculation by residues yields
 \begin{multline*}
 \label{9}
 \frac 1 {2i\pi}\langle f_k, g_{m,n} \rangle=\frac 1 {2i \pi}\int_{-\infty}^\infty
    \frac{h_k(t)e^{-2i\pi(n-k)t}}{t-(\alpha m{-}iw)}dt \\
               =  \begin{cases}
                   e^{-2i\pi n\alpha m} h_k(\alpha m-iw)e^{2i\pi k \alpha m}e^{2\pi w (k-n)} , & k < n; \\
                     0, & k\geq n,                  
              \end{cases}
\end{multline*}  
and
\beqq
\label{10} 
c_{m,n}=e^{-2i \pi n \alpha m}\sum_{k< n}h_k(\alpha m-iw)e^{2i\pi k \alpha m}e^{2\pi w(k-n)}. 
\eeqq        
 Let 
 \beq
 \label{11}
d_{m,n}:=c_{m,n}e^{2i \pi n \alpha m}, \ 
 \omega_{ m, k}:=h_k(\alpha m -iw)e^{2\pi i\alpha m k},
 \eeq
  and 
\beqq
{\vec d}_m=\{d_{m,n}\}_{n\in \ZZ},  \
 {\vec \omega}_m= \{\omega_{m, k}\}_{k\in \ZZ}.
 \eeqq

 \subsection*{c} {Norm estimates of the coefficients} 
We have
\beqq
\label{11c}
 \sum_m \|{\vec d}_m\|^2 = \|\{c_{m, n}\}\|^2_{l^2(\ZZ \times \ZZ)}, 
   \eeqq
 and 
\beq
\label{11}
A{\vec \omega}_m= {{\vec d}}_m, 
\eeq 
where the matrix $A$  is defined as
\beq
\label{11a}
A=(a_{n,k})_{n,k\in \ZZ}, \quad a_{n,k}= \begin{cases}
e^{2\pi(n-k)w}, & k < n; \\
0, & k\geq n. 
               \end{cases}.
\eeq

The matrix  $A$ defines the T{oe}plitz operator with the symbol
 \beqq
 s(\theta)=\sum_{k< 0} e^{2i \pi k \theta} e^{2\pi k w}= \frac {e^{-2i\pi \theta}}{e^{2\pi w} - e^{-2i \pi \theta}}, \ 
 \theta\in [0,2\pi].
 \eeqq
 We have 
 \beqq
 \max_{\theta \in [0,2\pi]} |s(\theta)|=\frac{1}{e^{2\pi w}-1}, \ 
 \min_{\theta \in [0,2\pi]} |s(\theta)|=\frac{1}{e^{2\pi w}+1}.
 \eeqq
 
Respectively
$$
(e^{2\pi w} -1) \|\vec {c}_m\|\leq\|\vec{\omega}_m\| \leq (e^{2\pi w} +1) \|\vec {c}_m\|,
$$
and
\beq
\label{eq1n}
 \left(e^{2\pi w} -1 \right)^2 \sum_m\|\vec {c}_m\|^2 \leq
\sum_m\|\vec{\omega}_m\|^2 \leq \left( e^{2\pi w} +1 \right)^2 \sum_m\|\vec {c}_m\|^2.
\eeq

 \subsection*{Frame constants} 
 
 On the other hand,
 \beq
\label{11c} 
 \sum_m \|{\vec \omega}_m\|^2=\sum_{  k}\sum_m |h_k(\alpha m-iw)|^2 , 
 \eeq
and it remains to compare the right-hand side to $\|f\|^2$.

Let $g_k(z)=h_k(z-iw)$. We have $g_k\in \PW\subset\PWa$. The Paley-Wiener representation yields  $$
\|h_k\|_{\Lt } \leq\|g_k\|_{\Lt}\leq  e^{2\pi w}\|h_k\|_{\Lt }.
$$

Finally
\beqq
\|g_k\|_{\Lt }^2= {\alpha} \sum |g_k(\alpha m)|^2,
\eeqq
this follows from the fact that the sequence of the reproducing kernels 
$$
\left \{  e^{i \frac \pi \alpha t  } \ \frac {\sin \frac \pi \alpha(t-\alpha m )}{\frac \pi \alpha(t-\alpha m)}
\right \}_{m\in \ZZ}
$$
is an {orthogonal} basis in $\PWa$.

By combining this with equations (\ref{5}), (\ref{eq1n}), (\ref{6}), and (\ref{7})
we obtain  
 \beqq
 \sum_{m,n} |c_{m,n}|^2 \geq \frac 1 {\alpha}\left (\frac{1}{1+e^{2\pi w}}\right )^2\|f\|^2,
 \eeqq
 and
  \beqq
 \sum_{m,n} |c_{m,n}|^2 \leq \frac 1 {\alpha}\left (\frac{e^{2\pi w}}{e^{2\pi w}-1}\right )^2\|f\|^2.
 \eeqq

Respectively  

\beqq
\label{low}
 A(\alpha,1;w)\geq   \frac {4\pi^2} {\alpha}\left (\frac{1}{1+e^{2\pi w}}\right )^2,
\eeqq
and 
\beqq
\label{upper}
 B(\alpha,1;w)\leq   \frac {4\pi^2} {\alpha} \left (\frac{e^{2\pi w}}{e^{2\pi w}-1}\right )^2.  
\eeqq

Now relation (\ref{norm}) yields Theorem \ref{main}.


\section{Frame operator and    the canonical dual window}\label{Section3}

 Recall that the frame operator $S=S^{\alpha,\beta}$ is defined as 
$$Sf=\sum_{m,n}\langle f,g_{m,n}\rangle g_{m,n}.$$
For $f\in C^\infty_0(\mathbb{R})$ we have 
$$Sf(x)=\sum_{m,n}\int_{\mathbb{R}}f(t)\frac{e^{2\pi i \beta m(x-t)}}{(t-\alpha n+iw)(x-\alpha n-iw)}dt.$$
Note that
$$\frac{1}{(t-\alpha n+iw)(x-\alpha n-iw)}=\frac{1}{x-t-2i w}\biggl{(}\frac{1}{t-\alpha n+iw}-\frac{1}{x-\alpha n-iw}\biggr{)},$$
and 
$$\sum_n\frac{1}{z_1-\alpha n}-\frac{1}{z_2-\alpha n}=\frac{\pi}{\alpha}\cot\biggl{(}\frac{\pi}{\alpha}z_1\biggr{)}-\frac{\pi}{\alpha}\cot\biggl{(}\frac{\pi}{\alpha}z_2\biggr{)}.$$
Hence,
$$Sf(x)=\frac{\pi}{\alpha}\int_\mathbb{R}f(t)\sum_m e^{2\pi i \beta m(x-t)}\frac{\cot\biggl{(}\frac{\pi}{\alpha}(t+i w)\biggr{)}-\cot\biggl{(}\frac{\pi}{\alpha}(x-iw)\biggr{)}}{x - t - 2iw}dt=$$
$$=\frac{\pi}{\alpha}\frac{1}{\sin\bigl{(}\frac{\pi}{\alpha}(x-iw)\bigr{)}}\int_\mathbb{R}\frac{f(t)}{\sin\bigl{(}\frac{\pi}{\alpha}(t+iw)\bigr{)}}\sum_m e^{2\pi i \beta m(x-t)}\frac{\sin(\frac{\pi}{\alpha}(x-t-2iw))}{x-t-2iw}dt.$$
  Consider (in the distributional sense) the 
 function $h$,
$$
h(s):=\sum_me^{2\pi i \beta m s}\frac{\sin(\frac{\pi}{\alpha}(s-2iw))}{s-2iw}.
$$
We have 
$$\hat{h}(\xi)={\pi}\sum_m\chi_{[-(2\alpha)^{-1},(2\alpha)^{-1}]+\beta m}(\xi)e^{4\pi w (\xi-\beta m)}.$$
In particular, $\hat{h}\simeq 1$ {since $\alpha \beta \le 1$.}
Therefore, 

$$Sf(x) =\frac{\pi}{\alpha}\sin^{-1}\left(\frac{\pi}{\alpha}(x-iw)\right)\left (f(t)\sin^{-1}\left(\frac{\pi}{\alpha}(t+iw)\right)\ast h\right).$$

Corrollary \ref{cor} immediately follows from this repesentation.

\subsection* {Canonical dual window}
We know that the canonical dual window $\gamma$ is given by $\gamma=S^{-1}g$.
We have 

$$S^{-1}f(t)=\frac{\alpha}{\pi}\sin\left(\frac{\pi}{\alpha}(t+iw)\right) \left(f(x)\sin\left(\frac{\pi}{\alpha}(x-iw)\right)\ast \tilde{h}\right),$$
where $\tilde{h}$ is such that $\hat{\tilde{h}}\hat{h}=1$.
  Then 
$$\gamma(t)=S^{-1}g(t)=\frac{\alpha}{\pi}\sin\left(\frac{\pi}{\alpha}(t+iw)\right) 
\left(\frac{\sin\left(\frac{\pi}{\alpha}(x-iw)\right)}{x-iw}\ast \tilde{h}\right)=$$
$$=
\frac{\alpha}{\pi}\sin\left(\frac{\pi}{\alpha}(t+iw)\right) \left(\frac{{\pi}\chi_{[-(2\alpha)^{-1},(2\alpha)^{-1}]}(\xi)e^{2\pi w \xi}}{\hat{h}(\xi)}\right)\check.$$
This gives us {a} formula for  $\gamma$. {Note that since the numerator is compactly supported, we can leave only finitely many terms in the formula for $\hat{h}(\xi)$. We will first do an explicit computation in the case $\alpha\beta \ge \frac{1}{2}$ and then explaining the changes 
 to be done in the
general case. For $\alpha\beta\geq \frac{1}{2}$ we have
$$\frac{{\pi}\chi_{[-(2\alpha)^{-1},(2\alpha)^{-1}]}(\xi)e^{2\pi w \xi}}{\hat{h}(\xi)}=$$
$$\frac{\chi_{[-(2\alpha)^{-1},(2\alpha)^{-1}]}(\xi)e^{2\pi w \xi}}{\sum_{m=-1,0,1}\chi_{[-(2\alpha)^{-1},(2\alpha)^{-1}]+\beta m}(\xi)e^{4\pi w (\xi-\beta m)}}=$$

$$ =\chi_{|\xi|<\beta-\frac{1}{2\alpha}}(\xi)e^{-2\pi w\xi}+\frac{1}{1+e^{-4\pi w\beta}}\chi_{[\beta-\frac{1}{2\alpha},\frac{1}{2\alpha}]}(\xi)e^{-2\pi w\xi}+\frac{1}{1+e^{4\pi w\beta}}\chi_{[\frac{-1}{2\alpha},-\beta+\frac{1}{2\alpha}]}(\xi)e^{-2\pi w\xi}.$$
{Collecting everything, we get}

$$\gamma(t)=\frac{\alpha}{\pi} \sin\left(\frac{\pi}{\alpha}(t+iw)\right)\biggl{[}
\frac{\sin(\pi(\beta-\varepsilon)(t{+}iw))}{{\pi(t+iw)}}+ {\left(\frac{e^{\pi w \beta - \pi i \beta t}}{1+e^{4\pi w \beta}}+\frac{e^{-\pi w \beta + \pi i \beta t}}{1+e^{-4\pi w \beta}}\right)}\frac{\sin(\pi\varepsilon(t{+}iw))}{{\pi(t+iw)}}\biggr{]},$$
where $\varepsilon=\frac{1}{\alpha}-\beta$.

{For the case $\frac{1}{n+1}\le \alpha\beta \le \frac{1}{n}$ { a} naive approach would give us a formula with $2n+1$ terms, and each term itself has a sum in the denominator. However, the sum in the denominator can be seen to be a geometric sum, which can be computed explicitly, and then summing odd and even terms on the Fourier transform side we get another geometric progression. Due to this answer also changes depending on whether $n$ is even or odd.}

We begin with the case of even $n = 2k$. Put $\lambda = \frac{1}{2\alpha} - k\beta$ and note that $0 \le \lambda \le \frac{\beta}{2}$. For $-\lambda < \xi < \lambda$ the formula for $\hat{h}(\xi)$ has $n+1$ non-zero terms, for $\lambda < \xi < \beta -\lambda$ it has $n$ non-zero terms and then it is periodic with the period $\beta$.

For $-\lambda < \xi < \lambda$ we have

$$\hat{h}(\xi) = \pi \sum_{m = -k}^k e^{4\pi w(\xi - \beta m)} =e^{4\pi w \xi}\pi e^{-4\pi k\beta w} \frac{e^{4\pi \beta w(2k+1)}-1}{e^{4\pi \beta w} - 1},$$

while for $\lambda < \xi < \beta-\lambda$ we have

$$\hat{h}(\xi) = \pi \sum_{m = -k+1}^k e^{4\pi w(\xi - \beta m)} = e^{4\pi w \xi}\pi e^{-4\pi (k-1)\beta w} \frac{e^{4\pi \beta w(2k) }-1}{e^{4\pi \beta w} - 1}.$$

We consider { the } intervals $[-\lambda + m\beta, \lambda + m\beta]$ for $-k\le m \le k$ and denote their union by $J$. We want to compute the inverse Fourier transform of $\frac{\pi \chi_J(\xi) e^{2\pi w\xi}}{\hat{h}(\xi)}$. We have

$$\frac{\pi \chi_J(\xi) e^{2\pi w\xi}}{\hat{h}(\xi)} =  e^{4\pi k\beta w} \frac{e^{4\pi \beta w} - 1}{e^{4\pi \beta w(2k+1)}-1}\sum_{m = -k}^k \chi_{[-\lambda + m\beta, \lambda + m\beta]}e^{2\pi w\xi}e^{-4\pi w(\xi - m\beta)}.$$

The first term is a constant so let us compute the inverse Fourier transform of the second term.

\begin{multline}\label{odd part}
\sum_{m = -k}^k \int_{-\lambda + m\beta}^{\lambda+m\beta} e^{2\pi i t\xi} e^{-2\pi w \xi} e^{4\pi w m \beta}d\xi = \int_{-\lambda}^\lambda e^{2\pi i t \xi - 2\pi w\xi} \sum_{m = -k}^k e^{2\pi it \beta m + 2\pi w\beta m}d\xi = \\
e^{-2\pi \beta k(w+it)}\frac{e^{2\pi \beta (2k+1)(w+it)}-1}{e^{2\pi \beta(w+it)}-1}\int_{-\lambda}^\lambda e^{(2\pi i t - 2\pi w)\xi}d\xi =\\ e^{-2\pi \beta k(w+it)}\frac{e^{2\pi \beta (2k+1)(w+it)}-1}{e^{2\pi \beta(w+it)}-1} \frac{\sin(2\pi \lambda(t+iw))}{\pi(t+iw)}.
\end{multline}

Now, we look at the intervals $[\lambda + m\beta, \beta - \lambda + m\beta]$ for $-k \le m \le k - 1$ and denote their union by $I$. Note that $I$ and $J$ are disjoint up to endpoints and together they are exactly covering the interval $[-\frac{1}{2\alpha}, \frac{1}{2\alpha}]$. Again, we compute the inverse Fourier transform of $\frac{\pi \chi_I(\xi) e^{2\pi w\xi}}{\hat{h}(\xi)}$. We have

$$\frac{\pi \chi_I(\xi) e^{2\pi w\xi}}{\hat{h}(\xi)} =   e^{4\pi (k-1)\beta w} \frac{e^{4\pi \beta w} - 1}{e^{4\pi \beta w(2k) }-1}\sum_{m = -k+1}^k \chi_{[\lambda+m\beta, \beta - \lambda + m\beta]}e^{2\pi w \xi}e^{-4\pi w(\xi - m\beta)}.$$

The first term is again { a} constant, so we compute the inverse Fourier transform of the second term.
\begin{multline}\label{even part}
\sum_{m = -k}^k \int_{\lambda + m\beta}^{\beta-\lambda+m\beta} e^{2\pi i t\xi} e^{-2\pi w \xi} e^{4\pi w m \beta}d\xi = \int_{\lambda}^{\beta-\lambda} e^{2\pi i t \xi - 2\pi w\xi} \sum_{m = -k+1}^k e^{2\pi it \beta m + 2\pi w\beta m}d\xi = \\
e^{-2\pi \beta (k-1)(w+it)}\frac{e^{2\pi \beta (2k)(w+it)}-1}{e^{2\pi \beta(w+it)}-1}\int_{\lambda}^{\beta-\lambda} e^{(2\pi i t - 2\pi w)\xi}d\xi =\\ e^{-2\pi \beta (k-1)(w+it)}\frac{e^{2\pi \beta (2k)(w+it)}-1}{e^{2\pi \beta(w+it)}-1} e^{(\pi i t - \pi w)\beta}\frac{\sin(\pi(\beta - 2\lambda)(t+iw))}{\pi(t+iw)}.
\end{multline}

Summing \eqref{odd part} and \eqref{even part} with  the corresponding constants, substituting $2k$ for $n$ everywhere, and then multiplying the result by $\frac{\alpha}{\pi}\sin\left(\frac{\pi}{\alpha}(t+iw)\right)$, we get
\begin{multline}\label{even n}
\gamma(t) = \Biggl (e^{2\pi (n-2)\beta w} \frac{e^{4\pi \beta w} - 1}{e^{4\pi \beta w n }-1} e^{-\pi \beta (n-2)(w+it)}\frac{e^{2\pi \beta n(w+it)}-1}{e^{2\pi \beta(w+it)}-1} e^{(\pi i t - \pi w)\beta}\frac{\sin(\pi(\beta - 2\lambda)(t+iw))}{\pi(t+iw)} + \\e^{2\pi n\beta w} \frac{e^{4\pi \beta w} - 1}{e^{4\pi \beta w(n+1)}-1}e^{-\pi \beta n(w+it)}\frac{e^{2\pi \beta (n+1)(w+it)}-1}{e^{2\pi \beta(w+it)}-1} \frac{\sin(2\pi \lambda(t+iw))}{\pi(t+iw)}\Biggr)\frac{\alpha}{\pi}\sin\left(\frac{\pi}{\alpha}(t+iw)\right) .
\end{multline}

Now we turn to the case of odd $n = 2k+1$. Put $\delta = (k+1)\beta - \frac{1}{2\alpha}$. Note that $0 \le \delta \le \frac{\beta}{2}$. For $-\delta < \xi < \delta$ we have

$$\hat{h}(\xi) = \pi \sum_{m = -k}^k e^{4\pi w(\xi - \beta m)} =e^{4\pi w \xi}\pi e^{-4\pi k\beta w} \frac{e^{4\pi \beta w(2k+1)}-1}{e^{4\pi \beta w} - 1},$$
while for $\delta < \xi < \beta - \delta$ we have

$$\hat{h}(\xi) = \pi \sum_{m = -k}^{k+1} e^{4\pi w(\xi - \beta m)} = e^{4\pi w \xi}\pi e^{-4\pi k\beta w} \frac{e^{4\pi \beta w(2k+2) }-1}{e^{4\pi \beta w} - 1}.$$

The intervals constituting $J$ now have the form $[-\delta + m\beta, \delta + m\beta]$ for $-k\le m \le k$, and the intervals constituting $I$ have the form $[\delta + m\beta, \beta - \delta + m\beta]$ for $-k-1 \le m \le k$ (note the increase by one from both sides in the latter case). The rest of the computation is exactly the same so we just present the final answer
\begin{multline}\label{odd n}
\gamma(t) = \Biggl( e^{2\pi (n-1)\beta w} \frac{e^{4\pi \beta w} - 1}{e^{4\pi \beta w(n+1)}-1} e^{-\pi \beta (n+1)(w+it)}\frac{e^{2\pi \beta (n+1)(w+it)}-1}{e^{2\pi \beta(w+it)}-1}e^{(\pi i t - \pi w)\beta}\frac{\sin(\pi(\beta - 2\delta)(t+iw))}{\pi(t+iw)} +\\ e^{2\pi (n-1)\beta w} \frac{e^{4\pi \beta w} - 1}{e^{4\pi \beta wn}-1}e^{-\pi \beta (n-1)(w+it)}\frac{e^{2\pi \beta n(w+it)}-1}{e^{2\pi \beta(w+it)}-1} \frac{\sin(2\pi \delta(t+iw))}{\pi(t+iw)}\Biggr)\frac{\alpha}{\pi}\sin\left(\frac{\pi}{\alpha}(t+iw)\right).
\end{multline}

\end{document}